\documentclass[12pt]{article} 
\usepackage{a4}
\usepackage{amssymb}
\usepackage{amsmath}
\newtheorem{theorem}{Theorem}[section]

\newtheorem{lemma}{Lemma}[section]

\newtheorem{definition}{Definition}[section]

\newtheorem{remark}{Remark}[section]
\begin{document}
\bibliographystyle{plain}
\title{Local existence and exponential growth for a semilinear damped wave equation 
with dynamic boundary conditions}
\author{St\'{e}phane Gerbi\thanks{Laboratoire de Math\'ematiques, Universit\'e de Savoie,  
73376 Le Bourget du Lac, France, e-mail:Stephane.Gerbi@univ-savoie.fr, \bf{corrresponding author.}} ~and 
Belkacem Said-Houari\thanks{ Laboratoire de Math\'{e}matiques Appliqu\'ees, Universit\'{e} Badji Mokhtar,
B.P. 12 Annaba 23000, Alg\'erie, e-mail:saidhouarib@yahoo.fr}}
\date{}
\maketitle

\begin{abstract}
In this paper we consider a multi-dimensional damped semiliear 
wave equation with dynamic boundary conditions, related to the Kelvin-Voigt damping. 
We firstly prove the local existence by using the Faedo-Galerkin approximations combined with a contraction mapping theorem. 
Secondly, the exponential growth of the energy and the $L^p$ norm of the solution is presented.
\end{abstract}

\section{Introduction}

In this paper we consider the following semilinear damped wave equation  with dynamic boundary conditions:
\begin{equation}\label{ondes}
\left\{
\begin{array}{ll}
u_{tt}-\Delta u-\alpha \Delta u_{t}=\vert u\vert^{p-2}u, & x\in \Omega ,\ t>0 \\[0.1cm]
u( x,t) =0, &  x\in \Gamma_{0},\ t>0  \\[0.1cm]
u_{tt}(x,t) =- \left[\displaystyle \frac{\partial u}{\partial \nu }(x,t) +
\frac{\alpha \partial u_{t}}{\partial \nu }(x,t) +r\vert u_{t}\vert^{m-2}u_{t}( x,t) \right]
&  x\in \Gamma_{1},\ t>0 \\[0.1cm]
u( x,0) =u_{0}(x), \; u_{t}( x,0) =u_{1}(x) & x\in \Omega \quad .
\end{array}
\right.
\end{equation}
where $u=u(x,t)\,,\, t \geq 0\,,\, x\in \Omega\,,\,\Delta$ denotes the Laplacian operator with res\-pect 
to the $x$ variable,  $\Omega$ is a regular and bounded domain of $\mathbb{R}^{N}\,,\,( N\geq 1)$,
$\partial\Omega~=~\Gamma _{0}~\cup~\Gamma _{1}$, $mes(\Gamma _{0}) >0,$ 
$\Gamma _{0}\cap \Gamma _{1}=\varnothing $ and $\displaystyle \frac{\partial }{\partial \nu }$ 
denotes the unit outer normal derivative, $m \geq 2\,,\,a\,,\,\alpha \mbox{ and } r$ are
positive constants, $p>2 $ and $u_{0}\,,\,u_{1}$ are given functions.

From the mathematical point of view, these problems do not neglect  acceleration terms on the boundary.
Such type of boundary conditions are usually called \textit{dynamic boundary conditions}.
They are not only important from the theoretical point of view but also arise in
several physical applications. In one space dimension, the problem (\ref{ondes}) can modelize the dynamic evolution of a viscoelastic rod that is fixed at one end and has a tip mass attached to its free end.
The dynamic boundary conditions represents  the Newton's law for the attached mass, 
(see \cite{BST64,AKS96, CM98} for more details). 
In the two dimension space, as showed in \cite{G06} and in the references therein, 
these  boundary conditions arise  when we consider the transverse motion of a flexible membrane $\Omega$ whose boundary 
may be affected by the vibrations only in a region. 
Also some dynamic boundary conditions as in problem (\ref{ondes}) appear when
we assume that $\Omega$ is an exterior domain of $\mathbb{R}^3$ in which homogeneous fluid is at rest except for sound waves.
Each point of the boundary is subjected to small normal displacements into the obstacle (see \cite{B76} for more details).
This type of dynamic boundary conditions are known as acoustic boundary conditions. 

In the one dimensional case and for $r=0,$ that is in the absence of
boundary damping, this problem has been considered by Grobbelaar-Van Dalsen \cite{G94_1}. 
By using the theory of $B$-evolutions and the theory of
fractional powers developped in \cite{S81,S97}, the author showed that the partial differential equations in the problem
(\ref{ondes})  gives rise to an analytic semigroup in
an appropriate functional space. As a consequence, the existence and the
uniqueness of solutions was obtained.
In the case where $r\neq 0 \mbox{ and } m = 2$, Pellicer and Sol{\`a}-Morales 
\cite{PS04} considered the one dimensional problem as an alternative model for the classical spring-mass damper system, 
and by using the dominant eigenvalues method, they proved that for small values of
the parameter $a$ the partial differential equations in the problem (\ref{ondes}) has the classical second order
differential equation
\[
m_1 u''(t)+ d_1 u'(t)+k_1 u(t)=0,
\]
as a limit where the parameter $m_1 \,,\, d_1 \mbox{ and } k_1$ are determined from the values of the spring-mass damper system. Thus, the asymptotic stability of the model has been determined as a consequence of this limit.
But they did not obtain any rate of convergence.

We recall that the presence of the strong damping term $-\Delta u_{t}$ in the problem (\ref{ondes}) makes the problem 
different from that considered in \cite{GT94} and widely studied in the litterature \cite{V99,T98,T99,GS06,TV05} for instance. 
For this reason less results were known for the wave equation with a strong damping and many problems remained unsolved,  specially the blow-up of solutions in the presence of a strong damping and nonlinear damping at the same time. 
Here we will give a partial answer to this question. 
That is to say, we will prove that the solution is unbounded and grows up exponentially when time goes to infinity.

Recently, Gazzola and Squassina \cite{GS06} studied the global solution and the finite time blow-up for
a damped semilinear wave equations with Dirichlet boundary conditions by a careful study 
of the stationnary solutions and their stability using the Nehari manifold and a mountain pass energy level of the initial
condition.

The main difficulty of the problem considered is related to the non ordinary boundary conditions defined on $\Gamma_1$. 
Very little attention has been paid to this type of boundary conditions. We mention only a few particular results in the one dimensional space and for a linear damping i.e. $(m=2)$ \cite{GV99,PS04,DL02}.

A related problem to (\ref{ondes}) is the following:
\begin{eqnarray*}
u_{tt}-\Delta u + g(u_{t}) &=&f \hspace*{1.5cm} \text{ in } \Omega \times ( 0,T)  \nonumber \\
\frac{\partial u}{\partial \nu }+ K(u) u_{tt}+ h(u_{t}) &=&0,\hspace*{1.5cm} 
\text{ on }\partial \Omega \times (0,T)\\
u(x,0) &=& u_{0}(x) \hspace*{1cm}\text{ in }\Omega \nonumber \\
u_{t}(x,0)&=&u_{1}(x) \hspace*{1cm}\text{ in }\Omega  \nonumber
\end{eqnarray*}
where the boundary term $h( u_{t}) =\left\vert u_{t}\right\vert
^{\rho }u_{t}$ arises when one studies flows of gaz in a channel with
porous walls. The term $u_{tt}$ on the boundary appears from the internal
forces, and the nonlinearity $K(u) u_{tt}$ on the boundary represents the internal
forces when the density of the medium depends on the displacement.
This problem has been studied in \cite{DL02,DLS98}. By
using the Fadeo-Galerkin approximations and a compactness argument they
proved the global existence and the exponential decay of the solution of
the problem.

We recall some results related to the interaction of an elastic medium with rigid mass. By using the classical
semigroup theory, Littman and Markus \cite{LM88} established a uniqueness result for a particular Euler-Bernoulli beam rigid
body structure. They also proved the asymptotic stability of the structure by using the feedback boundary
damping. In \cite{LL98} the authors considered the Euler-Bernoulli beam equation which describes the dynamics of clamped
elastic beam in which one segment of the beam is made with viscoelastic material and the other of elastic
material. By combining the frequency domain method with the multiplier technique, they proved the exponential decay 
for the transversal motion but not for the longitudinal motion of the model, when the Kelvin-Voigt damping
is distributed only on a subinterval of the domain. 
In relation with this point, see also the work by Chen et \textit{al.} \cite{CLL98} concerning the Euler-Bernoulli beam
equation with the global or local Kelvin-Voigt damping. Also models of vibrating strings with local viscoelasticity and
Boltzmann damping, instead of the Kelvin-Voigt one, were considered in \cite{LL02} and an exponential energy decay 
rate was established. Recently, Grobbelaar-Van Dalsen \cite{G03} considered an extensible thermo-elastic beam 
which is hanged at one end with rigid body attached to its free end, i.e. one dimensional hybrid thermoelastic structure,
and showed that the method used in \cite{O97} is still valid to establish an uniform stabilization of the system. 
Concerning the controllability of the hybrid system we refer to the work by Castro and Zuazua\cite{CZ98}, in which they
considered flexible beams connected by point mass and the model takes account of the rotational inertia.

In this paper we consider the problem (\ref{ondes}) where we have set
for the sake of simplity $a = 1$. Section 2 is devoted to the local existence and uniqueness of the
solution of the problem (\ref{ondes}). We will use a technique close to the one used by Georgiev and Todorova in
\cite{GT94} and Vitillaro in \cite{V02_1,V02_2}: a Faedo-Galerkin approximation coupled to a fix point theorem.

In section 3, we shall prove that the energy is unbounded  when the initial data are large enough. 
In fact, it will be proved that the $L^{p}$-norm of the solutions grows as an exponential
function. An essential ingredient of the proof is a lower bound in the $L^p$ norm and the $H^1$ 
seminorm of the solution when the initial data are large enough, obtained by  Vitillaro in \cite{V99}.
The other ingredient is the use of an auxillary function $L$ (which is a small  perturbation of the energy)
in order to obtain a linear differential inequality, that we integrate to finally prove that the energy is 
exponentially growing. To this end, we use Young's inequality with suitable coefficient, interpolation 
and Poincar\'e's inequalities.
 
Let us recall that the blow-up result in the case of a nonlinear
damping $( m \neq 2) $ is still an open problem.
\section{Local existence}
In this section we will prove the local existence and the uniqueness of the solution of the problem (\ref{ondes}).
We will adapt the ideas used by Georgiev and Todorova in \cite{GT94}, which consists in constructing approximations by the
Faedo-Galerkin procedure in order to use the contraction mapping theorem.
This method allows us to consider less restrictions on the initial data. Consequently, the same result can be established by using the Faedo-Galerkin approximation method coupled with the potential well method \cite{CCM04}.
\subsection{Setup and notations}
We present here some material that we shall use in order to prove the local existence of the solution of
problem (\ref{ondes}).
We denote
$$
H_{\Gamma_{0}}^{1}(\Omega) =\left\{u \in H^1(\Omega) /\ u_{\Gamma_{0}} = 0\right\} .
$$
By $( .,.) $ we denote the scalar product in $L^{2}( \Omega)$  i.e. $(u,v)(t) = \displaystyle \int_{\Omega} u(x,t) v(x,t) dx$. Also we mean by $\Vert .\Vert_{q}$ the $L^{q}(\Omega) $ norm for $1 \leq q \leq \infty$, and by 
$\Vert .\Vert_{q,\Gamma_{1}}$ the $L^{q}(\Gamma_{1}) $ norm.

Let $T>0$ be a real number and $X$ a Banach space endowed with norm $\Vert .\Vert _{X}$. 
$L^{p}(0,T;X) ,\ 1 \leq p < \infty$ denotes the space of functions $f$ which are $L^{p}$ 
over $\left( 0,T\right) $ with values in $X$, which are measurable and  $\Vert f \Vert_{X} \in L^{p} \left(0,T\right)$. 
This space is a Banach space endowed with the norm
$$
\Vert f\Vert_{L^{p}\left( 0,T;X\right) }=
\left(\int_{0}^{T}\Vert f\Vert_{X}^{p} dt\right)^{1/p} \quad .
$$
$L^{\infty}\left( 0,T;X\right) $ denotes the space of functions $f:\left]0,T\right[ \rightarrow X$ which are measurable and 
$\Vert f\Vert_{X}\in L^{\infty }\left( 0,T\right) $.
This space is a Banach space endowed with the norm:
$$
\Vert f\Vert_{L^{\infty}(0,T;X)}=\mbox{ess}\sup_{0<t<T}\Vert f\Vert_{X} \quad .
$$
We recall that if $X$ and $Y$ are two Banach spaces such that $X\hookrightarrow Y$ (continuous embedding), then
$$
L^{p}\left( 0,T;X\right) \hookrightarrow L^{p}\left( 0,T;Y\right) , \ 1 \leq p\leq \infty .
$$
We will also use the embedding (see \cite[Therorem 5.8]{A75}).
$$
H_{\Gamma_{0}}^{1}(\Omega) \hookrightarrow L^{q}(\Gamma_{1}), \;  2 \leq q \leq \bar{q} \quad
\mbox{where } \quad 
\bar{q}=\left\{
\begin{array}{c}
\displaystyle \frac{2 (N-1) }{N-2} \,,\mbox{ if } N \geq 3 \\
+\infty \,,\mbox{ if } N=1, 2
\end{array}
\right. .
$$
Let us denote $V = H_{\Gamma_{0}}^{1}(\Omega) \cap L^{m}(\Gamma_{1}) $.

In this work, we cannot use ``directly'' the existence result of Georgiev and Todorova \cite{GT94} nor the results of Vitillaro
\cite{V02_1,V02_2} because of the presence of the strong linear damping $- \Delta u_t$ and the dynamic boundary
conditions on $\Gamma_1$.
Therefore, we have the next local existence theorem.
\begin{theorem}\label{existence}
Let $2\leq p\leq \bar{q}$ and $\max\left( 2,\displaystyle \frac{\bar{q}}{\bar{q}+1-p} \right) \leq m \leq \bar{q}$. \\
Then given $u_{0}\in H_{\Gamma_{0}}^{1}(\Omega) $ and $u_{1}\in L^{2}(\Omega) $, there exists $T > 0$ and a unique
solution $u$ of the problem (\ref{ondes}) on $(0,T) $ such that
\begin{eqnarray*}
u &\in &C \Bigl( [ 0,T],H_{\Gamma_{0}}^{1}(\Omega) \Bigl) \cap C^{1}\Bigl( [ 0,T] ,L^{2}(\Omega) \Bigl), \\
u_{t} &\in &L^{2}\Bigl(0,T;H_{\Gamma_{0}}^{1}(\Omega)\Bigl) \cap L^{m}\Bigl( (0,T) \times \Gamma_{1}\Bigl)
\end{eqnarray*}
\end{theorem}
We will prove this theorem by using the Fadeo-Galerkin approximations and the well-known contraction mapping theorem. 
In order to define the function for which a fixed point exists, we will consider first a related problem. \\
For $u  \in C \bigl( [ 0,T],H_{\Gamma_{0}}^{1}(\Omega) \bigl) \, \cap \, C^{1}\bigl( [ 0,T] ,L^{2}(\Omega) \bigl)$ given, 
let us consider the following problem:
\begin{equation}\label{ondes_u}
\left\{
\begin{array}{ll}
v_{tt}-\Delta v-\alpha \Delta v_{t}=\vert u\vert^{p-2}u, & x\in \Omega ,\ t>0 \\[0.1cm]
v( x,t) =0, &  x\in \Gamma_{0},\ t>0  \\[0.1cm]
v_{tt}(x,t) =- \left[\displaystyle \frac{\partial v}{\partial \nu }(x,t) +
\frac{\alpha \partial v_{t}}{\partial \nu }(x,t) +r\vert v_{t}\vert^{m-2}v_{t}( x,t) \right]
&  x\in \Gamma_{1},\ t>0 \\[0.1cm]
v(x,0) =u_{0}(x), \; v_{t}( x,0) =u_{1}(x) & x\in \Omega \quad .
\end{array}
\right.
\end{equation}
We have now to state the following existence result:
\begin{lemma}\label{existence_u}
Let $2\leq p\leq \bar{q}$ and $\max\left( 2,\displaystyle \frac{\bar{q}}{\bar{q}+1-p} \right) \leq m \leq \bar{q}$. 
Then given $u_{0}~\in~H_{\Gamma_{0}}^{1}(\Omega)\,,\, u_{1}\in L^{2}(\Omega)$ there exists $T > 0$ and a 
unique solution $v$ of the problem (\ref{ondes_u}) on $(0,T) $ such that
\begin{eqnarray*}
v &\in &C \Bigl( [ 0,T],H_{\Gamma_{0}}^{1}(\Omega) \Bigl) \, \cap \,C^{1}\Bigl( [ 0,T] ,L^{2}(\Omega) \Bigl), \\
v_{t} &\in &L^{2}\Bigl(0,T;H_{\Gamma_{0}}^{1}(\Omega)\Bigl) \,\cap \,L^{m}\Bigl( (0,T) \times \Gamma_{1}\Bigl)
\end{eqnarray*}
and satisfies the energy identity:
\begin{eqnarray*}
\frac{1}{2}\left[\Vert \nabla v\Vert_{2}^{2}+\Vert v_{t}\Vert_{2}^{2}+\Vert v_{t}\Vert_{2,\Gamma_{1}}^{2}\right]_{s}^{t}& + &
\alpha \displaystyle \int_{s}^{t}\Vert \nabla v_{t}(\tau)\Vert_{2}^{2} d\tau + 
r \displaystyle \int_{s}^{t}\Vert v_{t}(\tau)\Vert_{m,\Gamma_{1}}^{m} d\tau \\
&=&\displaystyle \int_{s}^{t} \displaystyle \int_{\Omega} \vert u(\tau)\vert^{p-2} u(\tau)v_{t}(\tau) d\tau dx
\end{eqnarray*}
for $0\leq s\leq t\leq T$.
\end{lemma}
In order to prove lemma \ref{existence_u}, we first study for any $T>0$ and $f\in H^{1}( 0,T;L^{2}(\Omega))$
the following problem:
\begin{equation}\label{ondes_f}
\left\{
\begin{array}{ll}
v_{tt}-\Delta v-\alpha \Delta v_{t}=f(x,t), & x\in \Omega ,\ t>0 \\[0.1cm]
v( x,t) =0, &  x\in \Gamma_{0},\ t>0  \\[0.1cm]
v_{tt}(x,t) =- \left[\displaystyle \frac{\partial v}{\partial \nu }(x,t) +
\frac{\alpha \partial v_{t}}{\partial \nu }(x,t) +r\vert v_{t}\vert^{m-2}v_{t}( x,t) \right]
&  x\in \Gamma_{1},\ t>0 \\[0.1cm]
v(x,0) =u_{0}(x), \; v_{t}( x,0) =u_{1}(x) & x\in \Omega \quad .
\end{array}
\right.
\end{equation}
At this point, as done by Doronin et \textit{al.} \cite{DLS98}, we have to precise exactly what type of solutions 
of the problem (\ref{ondes_f}) we expected. 
\begin{definition} \label{generalised}
A function $v(x,t) $ such that
\begin{eqnarray*}
v &\in &L^{\infty}\left( 0,T;H_{\Gamma_{0}}^{1}(\Omega) \right) \ , \\
v_{t} &\in &L^{2}\left( 0,T;H_{\Gamma_{0}}^{1}(\Omega) \right) \cap L^{m}\left( ( 0,T) \times \Gamma_{1}\right) \ , \\
v_{t} &\in &L^{\infty}\left( 0,T;H_{\Gamma_{0}}^{1}(\Omega)\right) \cap L^{\infty}\left( 0,T;L^{2}(\Gamma_{1})\right)\ , \\
v_{tt} &\in &L^{\infty}\left( 0,T;L^{2}(\Omega)\right) \cap L^{\infty}\left( 0,T;L^{2}(\Gamma_{1}) \right) \ , \\
v(x,0) & = &u_{0}(x)\,,\\ 
v_{t}(x,0) &= & u_{1}(x) \,,
\end{eqnarray*}%
is a generalized solution to the problem (\ref{ondes_f}) if for any function 
$\omega \in H_{\Gamma_{0}}^{1}(\Omega) \cap L^{m}(\Gamma_{1}) $ and $\varphi \in C^{1}(0,T) $ with 
$\varphi(T) =0$, we have the following identity:
\begin{equation*}
\begin{array}{lll}
\hspace*{-0.2cm}\displaystyle  \int_{0}^{T} (f,w)(t) \, \varphi(t)\, dt &=& \displaystyle \int_{0}^{T}
\Bigl[ (v_{tt},w)(t) +(\nabla v,\nabla w)(t) + \alpha ( \nabla v_{t},\nabla w)(t)\Bigl] \, \varphi(t)\, dt \\
&+&\displaystyle \int_{0}^{T}\varphi(t) \int_{\Gamma_{1}}
\Bigl[v_{tt}(t) +r \vert v_{t}(t) \vert^{m-2}v_{t}(t) \Bigl] w \, d\sigma \, dt .
\end{array}%
\end{equation*}
\end{definition}
\begin{lemma}\label{existence_f}
Let $2\leq p\leq \bar{q}$ and $2 \leq m \leq \bar{q}$. \\ 
Let $u_{0}\in H^2(\Omega) \cap V \,,\, u_{1}\in H^{2}(\Omega)$ and $f\in H^{1}( 0,T;L^{2}(\Omega))$, then 
for any $T > 0,$ there exists a unique generalized solution (in the sense of definition \ref{generalised}),
$v(t,x) $ of problem (\ref{ondes_f}).
\end{lemma}
\subsection{Proof of the lemma \ref{existence_f}}
To prove the above lemma, we will use the Faedo-Galerkin method, which consists in constructing approximations of the 
solution, then we obtain a priori estimates necessary to guarantee the convergence of these approximations.
It appears some difficulties in order to derive a second order estimate of $v_{tt}(0)$. To get rid of them, and 
inspired by the ideas of Doronin and Larkin in \cite{DL02} and Cavalcanti et \textit{al.}
\cite{CCSM00}, we introduce the following change of variables:
$$
\widetilde{v}(t,x) =v(t,x) - \phi (t,x)
\mbox{ with }
\phi(t,x) = u_{0}(x) +t \, u_{1}(x) .
$$
Consequently, we have the following problem with the unknown $\widetilde{v}(t,x)$ and null initial conditions:
\begin{equation}\label{ondes_ftilde}
\left\{
\begin{array}{lll}
\widetilde{v}_{tt}-\Delta \widetilde{v} -\alpha \Delta \widetilde{v}_{t}&=f(x,t)
+ \Delta \phi + \alpha \Delta \phi_t, & x\in \Omega ,\ t>0 \\[0.1cm]
\widetilde{v}(x,t) = 0, & &  x\in \Gamma_{0},\ t>0  \\[0.1cm]
\widetilde{v}_{tt}(x,t) =- \Bigl[&\hspace*{-0.5cm} \displaystyle \frac{\partial(\widetilde{v} + \phi)}{\partial \nu }(x,t) +
\frac{\alpha \partial (\widetilde{v}_{t}+\phi_t)}{\partial \nu }(x,t) \Bigl] - & ~ \\
~&\Bigl( r\vert (\widetilde{v}_{t}+ \phi_t)\vert^{m-2}(\widetilde{v}_{t}+\phi_t)(x,t)\Bigl) 
&  x\in \Gamma_{1},\ t>0 \\[0.1cm]
\widetilde{v}(x,0) = 0, \;  &\widetilde{v}_{t}(x,0) = 0 & x\in \Omega \quad .
\end{array}
\right.
\end{equation}
\begin{remark}\rm 
It is quite clear that if $\widetilde{v}$ is a solution of problem (\ref{ondes_ftilde}) on $[0,T]$, then $v$ is a solution of
problem (\ref{ondes_f}) on $[0,T]$. Moreover writing the problem in term of $\widetilde{v}$ shows exactly the regularity needed
on the initial conditions $u_0 \text{ and } u_1$ to ensure the existence.
\end{remark}
Now we construct approximations of the solution $\widetilde{v}$ by the Faedo-Galerkin method as follows.\\
For every $n \geq 1$, let $W_{n}= \mbox{span}\{\omega_{1},\ldots,\omega_{n}\},$ where 
$\{\omega_{j}(x)\}_{1\leq j \leq n}$ is a basis in the space $V$.
By using the Grahm-Schmidt orthogonalization process we can take $\omega =\bigl(\omega_{1},\ldots,\omega_{n}\bigl) $
to be orthonormal\footnote{Unfortunately, the presence of the nonlinear boundary conditions excludes us to use the spatial basis
of eigenfunctions of $-\Delta$ in $H^1_{\Gamma_0}(\Omega)$ as done in \cite{GS06}} in 
$L^{2}(\Omega)\cap L^{2}(\Gamma_1)$.\\
We define the approximations:
\begin{equation} \label{approx}
\widetilde{v}_{n}(t) =\sum_{j=1}^{n} g_{jn}(t) w_{j}
\end{equation}
where $\widetilde{v}_{n}(t) $ are solutions to the finite dimensional Cauchy problem 
(written in normal form since $\omega$ is an orthonormal basis):
\begin{equation}\label{approxsyst}
\begin{array}{c}
\displaystyle \int_{\Omega}\widetilde{v}_{ttn}(t) w_{j} \, dx +
\displaystyle \int_{\Omega}\nabla\bigl(\widetilde{v}_{n}+\phi \bigl) \nabla w_{j} +
\alpha \displaystyle \int_{\Omega}\nabla \bigl(\widetilde{v}_{n}+\phi \bigl)_{t} \nabla w_{j} \, dx
\\
+\displaystyle \int_{\Gamma_{1}} \left( \widetilde{v}_{ttn}(t)+r \vert(\widetilde{v}_{n} + \phi)_{t}\vert^{m-2}
(\widetilde{v}_{n} + \phi)_{t}\right) w_{j} \, d\sigma
= \int_{\Omega}fw_{j} \, dx \quad .\\
g_{jn}(0) = g'_{jn}(0) =0,\ j=1,\ldots, n
\end{array}
\end{equation}
According to the Caratheodory theorem, see \cite{CL55}, the problem (\ref{approxsyst}) has solution 
$(g_{jn}(t))_{j=1,n} \in H^{3}(0,t_{n})$ defined on $[0,t_{n})$.
We need now to show:
\begin{itemize}
\item firstly that forall $ n \in \mathbb{N} \,,_, t_{n}=T$,
\item secondly that these approximations converge to a solution
of the problem (\ref{ondes_ftilde}).
\end{itemize}
To do this we need the two following a priori estimates: first-order a priori estimates
to prove the first point. But we will show that the presence of the nonlinear term $\vert u_t\vert^{m-2} u_t$ 
forces us to derive a second order a priori estimate to pass to the limit in the nonlinear term. Indeed the key tool 
in our proof is the Aubin-Lions lemma which uses the compactness of the embedding 
$H^{\frac{1}{2}}(\Gamma_{1})\hookrightarrow L^{2}(\Gamma_{1})$.

\subsubsection{First order a priori estimates}
Multiplying equation (\ref{approxsyst}) by $g'_{jn}(t)$, integrating over $(0,t) \times \Omega $ and using integration
by parts we get: for every $n\geq 1$, 
\begin{equation}\label{energyn} 
\begin{array}{lll}
&\displaystyle \frac{1}{2}&\left[ \Vert\nabla \widetilde{v}_{n}(t)\Vert_{2}^{2}+ \Vert \widetilde{v}_{tn}(t)\Vert_{2}^{2}+
\Vert \widetilde{v}_{tn}\Vert_{2,\Gamma_{1}}^{2}\right] +
\displaystyle \int\limits_{0}^{t}\displaystyle \int\limits_{\Omega}\nabla \phi \nabla \widetilde{v}_{n} \, dx \\ 
&+~\alpha& \displaystyle \int\limits_{0}^{t}\displaystyle \int\limits_{\Omega}\nabla \phi_{t}
\nabla \widetilde{v}_{tn} \, dx + 
\alpha \displaystyle \int\limits_{0}^{t} \Vert \nabla \widetilde{v}_{tn}(s)\Vert_{2}^{2} \, ds  \\
&+~r&\displaystyle \int\limits_{0}^{t}\displaystyle \int\limits_{\Gamma_{1}}\vert \left(\widetilde{v}_{n}+\phi \right)_{t}
\vert^{m-2}(\widetilde{v}_{n}+\phi)_{t}\widetilde{v}_{tn} \,d\sigma ds \\
 &&~ = \displaystyle \int\limits_{0}^{t}\displaystyle \int\limits_{\Omega}f(t,x) \widetilde{v}_{tn}(s) \,dx\,ds
\end{array}
\end{equation}
By using Young's inequality, there exists $\delta_{1}>0$, (in fact small enough) such that
\begin{equation}\label{lhs1}
\alpha \displaystyle \int\limits_{0}^{t}\displaystyle \int\limits_{\Omega}\nabla \phi_{t}\nabla \widetilde{v_{t}}_{n}dx
\leq \delta_{1}\displaystyle \int\limits_{0}^{t}\displaystyle \int\limits_{\Omega}
\vert \nabla \widetilde{v_{t}}_{n}\vert^{2}dx + \frac{1}{4 \delta_{1}}
\displaystyle \int\limits_{0}^{t}\displaystyle \int\limits_{\Omega}\vert \nabla \phi_{t}\vert^{2}dx
\end{equation}
and
\begin{equation}\label{lhs2}
\displaystyle \int\limits_{0}^{t}\displaystyle \int\limits_{\Omega}\nabla \phi \nabla \widetilde{v}_{n}dx
\leq \delta_{1}\displaystyle \int\limits_{0}^{t}\displaystyle \int\limits_{\Omega}\vert\nabla \widetilde{v}_{n}\vert^{2} dx
+\frac{1}{4\delta_{1}}\displaystyle \int\limits_{0}^{t}\displaystyle \int\limits_{\Omega}\nabla \phi\vert^{2}dx.
\end{equation}
By Young's and Poincar\'{e}'s inequalities, we can find $C>0,$ such that
\begin{equation}\label{lhs3}
\displaystyle \int\limits_{0}^{t}\displaystyle \int\limits_{\Omega}f(t,x) \widetilde{v}_{tn}(s) dx ds\leq
C\displaystyle \int\limits_{0}^{t}\displaystyle \int\limits_{\Omega}\left(f^{2}+\vert \nabla \widetilde{v}_{tn}(s) \vert^{2}\right) dx ds.
\end{equation}
The last term in the left hand side of equation (\ref{energyn}) can be written as follows:
\begin{eqnarray*}
&&\displaystyle \int\limits_{0}^{t}\displaystyle \int\limits_{\Gamma_{1}}
\vert \left(\widetilde{v}_{n}+\phi \right)_{t}\vert^{m-2}\left( \widetilde{v}_{n}+\phi\right)_{t}\widetilde{v}_{tn}d\sigma ds \\
&=&\displaystyle \int\limits_{0}^{t}\displaystyle 
\int\limits_{\Gamma_{1}}\vert \left(\widetilde{v}_{n}+\phi \right)_{t}\vert^{m}d\sigma ds - 
\displaystyle \int\limits_{0}^{t}\displaystyle \int\limits_{\Gamma_{1}}
\vert \left(\widetilde{v}_{n}+\phi \right)_{t}\vert^{m-2}\left( \widetilde{v}_{n}+\phi \right)_{t}\phi_{t}d\sigma ds,
\end{eqnarray*}
Then Young's inequality gives us, for $\delta_{2}>0$:
\begin{equation}\label{lhs4}
\begin{array}{lll}
&&\left\vert \displaystyle \int\limits_{0}^{t}\displaystyle \int\limits_{\Gamma_{1}}
\vert \left( \widetilde{v}_{n}+\phi \right)_{t}\vert^{m-2}\left( \widetilde{v}_{n}+\phi\right)_{t}\phi_{t}d\sigma ds 
\right \vert \\
&\leq &\displaystyle \frac{\delta_{2}^{m}}{m}\int\limits_{0}^{t}\displaystyle \int\limits_{\Gamma_{1}}
\vert \left( \widetilde{v}_{n}+\phi \right)_{t}\vert^{m} d\sigma ds+
\frac{m-1}{m}\delta_{2}^{-m/\left( m-1\right)}\displaystyle \int\limits_{0}^{t}\displaystyle \int\limits_{\Gamma_{1}}
\vert \phi_{t}\vert^{m}d\sigma ds.  
\end{array}
\end{equation}
Consequently, using the inequalities (\ref{lhs1}), (\ref{lhs2}), (\ref{lhs3}) and (\ref{lhs4}) in the equation (\ref{energyn}),
choosing $\delta_{1}$ and $\delta_{2}$ small enough, we may conclude that:
\begin{equation}\label{1stestimate}
\begin{array}{lll}
&&\displaystyle \frac{1}{2}\left[ \Vert \nabla \widetilde{v}_{n}(t)\Vert_{2}^{2}+
\Vert \widetilde{v}_{tn}(t)\Vert_{2}^{2}+
\Vert \widetilde{v}_{tn}(t)\Vert_{2,\Gamma_{1}}^{2}\right]  \\
&&+\alpha \displaystyle \int\limits_{0}^{t}\Vert \nabla \widetilde{v}_{tn}(s) \Vert_{2}^{2}ds+
r\displaystyle \int\limits_{0}^{t}\displaystyle \int\limits_{\Gamma_{1}}\vert \left( \widetilde{v}_{n}+
\phi \right)_{t}\vert^{m}d\sigma ds\leq C_{T} \quad ,
\end{array}
\end{equation}
where $C_{T}$ is a positive constant independent of $n$.
Therefore, the last estimate (\ref{1stestimate}) gives us, $\forall n \in \mathbb{N} \,,\, t_n = T$, and:
\begin{equation}\label{regul1}
\left(\widetilde{v}_{n}\right)_{n\in \mathbb{N}} \text{ is bounded in }
L^{\infty}(0,T;H_{\Gamma_{0}}^{1}(\Omega)) \quad,
\end{equation}
\begin{equation}\label{regul1bis}
\begin{array}{ll}
\left(\widetilde{v}_{tn}\right)_{n\in \mathbb{N}} \text{ is bounded in }
L^{\infty}(0,T;L^{2}(\Omega)) \hspace*{-0.3cm}&\cap L^{2}\left(0,T;H_{\Gamma_{0}}^{1}(\Omega)\right) \\
\hspace*{-0.3cm}&\cap L^{\infty}\left(0,T;L^{2}(\Gamma_{1})\right) .
\end{array}
\end{equation}
Now, by using the following algebraic inequality:
\begin{equation}\label{algebr_lambda}
( A + B)^{\lambda }\leq 2^{\lambda -1}( A^{\lambda}+B^{\lambda }) ,\ A,B\geq 0\,,\, \lambda \geq 1,
\end{equation}
we can find $c_{1},c_{2}>0$, such that:
\begin{equation}\label{last1}
\displaystyle \int\limits_{0}^{t}\displaystyle \int\limits_{\Gamma_{1}}\vert (\widetilde{v}_{n}+\phi)_{t}\vert^{m}d\sigma ds
\geq c_{1}\displaystyle \int\limits_{0}^{t}\displaystyle \int\limits_{\Gamma_{1}}
\vert \widetilde{v}_{tn}\vert^{m}d\sigma ds - c_{2}\displaystyle \int\limits_{0}^{t}\displaystyle \int\limits_{\Gamma_{1}}
\vert \phi_{t}\vert^{m}d\sigma ds.
\end{equation}
Then, by the embedding $H^{2}(\Omega) \hookrightarrow L^{m}(\Gamma_{1})$ ($2 \leq m \leq \bar{q}$), 
we conclude that $u_{1}\in L^{m}(\Gamma_{1}) $.
Therefore, from the inequalities (\ref{1stestimate})  and (\ref{last1}), there exists $C'_{T}> 0 $ such that:
$$
\displaystyle \int\limits_{0}^{t}\displaystyle \int\limits_{\Gamma_{1}}\vert \widetilde{v}_{tn}\vert^{m}d\sigma ds
\leq C'_{T}.
$$
Consequently,
\begin{equation}\label{regul2}
\widetilde{v}_{tn} \text{ is bounded in } L^{m}((0,T) \times \Gamma_{1}) .
\end{equation}
\subsubsection{Second order a priori estimate}
In order to obtain a second a priori estimate, we will first estimate $\Vert \widetilde{v}_{ttn}(0) \Vert_{2}^{2}$ 
and $\Vert \widetilde{v}_{ttn}(0) \Vert_{2,\Gamma_{1}}^{2}$.
For this purpose, considering $w_{j}=\widetilde{v}_{ttn}(0) $ and $t=0$ in the equation (\ref{approxsyst}), we get
\begin{equation}\label{secondstep1}
\begin{array}{lll} 
&&\Vert \widetilde{v}_{ttn}(0) \Vert_{2}^{2}+\Vert \widetilde{v}_{ttn}(0) \Vert_{2,\Gamma_1}^{2}+
\displaystyle \int\limits_{\Omega}\nabla \phi(0) \nabla \widetilde{v}_{ttn}(0) dx   \\
&&+\alpha \displaystyle \int\limits_{\Omega}\nabla \phi_{t}(0)\nabla \widetilde{v}_{ttn}(0) dx+
r\displaystyle \int\limits_{\Gamma_{1}}\vert \phi_{t}(0)\vert^{m-2}\phi_{t}(0)\widetilde{v}_{ttn}(0) d\sigma ds  \\
&=&\displaystyle \int\limits_{\Omega}f(0,x) \widetilde{v}_{ttn}(0) dx ds.
\end{array}
\end{equation}
Since the following equalities hold:
$$
\phi(0) = u_0 \,,\, \phi_t(0) = u_1 \,,\, 
\displaystyle \int_\Omega \nabla \phi_t(0) \nabla \widetilde{v}_{ttn}(0) =
- \displaystyle \int_\Omega \Delta \phi_t(0) \widetilde{v}_{ttn}(0) + \int_{\Gamma_1} \phi_t \frac{\partial v_{ttn}}{\partial \nu} d\sigma \,,
$$
as $f\in H^{1}\left(0,T;L^{2}(\Omega)\right)$ and $u_{0},u_{1}\in H^{2}\left(\Omega \right)$, by using Young's inequality
and the embedding $H^{2}(\Omega) \hookrightarrow L^{m}(\Gamma_{1})$, we conclude that there exists $C>0$ 
independent of $n$ such that:
\begin{equation}\label{estim_0}
\Vert \widetilde{v}_{ttn}(0) \Vert_{2}^{2}+\Vert \widetilde{v}_{ttn}(0) \Vert_{2,\Gamma_{1}}^{2}\leq C \ .
\end{equation}
Differentiating equation (\ref{approxsyst}) with respect to $t$, multiplying the result by $g''_{jn}(t)$ 
and summing over $j$ we get:
\begin{equation}\label{diffsyst}
\begin{array}{lll}
&&\displaystyle \frac{1}{2}\frac{d}{dt}\left[ \Vert \nabla \widetilde{v}_{tn}(t) \Vert_{2}^{2}+
\Vert \widetilde{v}_{ttn}(t)\Vert_{2}^{2} + \Vert \widetilde{v}_{ttn}(t)\Vert_{2,\Gamma_{1}}^{2}\right] +
\displaystyle \int\limits_{\Omega}\nabla \phi_{t}\nabla \widetilde{v}_{ttn}dx \\
&&+\alpha \Vert \nabla \widetilde{v}_{ttn}(s) \Vert_{2}^{2}+r(m-1) \displaystyle \int\limits_{\Gamma_{1}}
\vert (\widetilde{v}_{n}+\phi)_{t}\vert^{m-2}(\widetilde{v}_{n}+\phi) _{tt}\widetilde{v}_{ttn}d\sigma\\
&=&\displaystyle \int\limits_{\Omega}f_{t}(t,x) \widetilde{v}_{ttn}(s) dx ds.
\end{array}
\end{equation}
Since $\phi_{tt}=0$, the last term in the left hand side of the equation (\ref{diffsyst}) can be written as follows:
$$
\displaystyle \int\limits_{\Gamma_{1}}
\vert (\widetilde{v}_{n}+\phi)_{t}\vert^{m-2}(\widetilde{v}_{n}+\phi) _{tt}\widetilde{v}_{ttn}d\sigma = 
\displaystyle \int\limits_{\Gamma_{1}}\vert( \widetilde{v}_{n}+\phi)_{t}\vert^{m-2}
(\widetilde{v}_{ttn}+\phi_{tt}) ^{2}d\sigma ,
$$
But we have,
$$
\displaystyle \int\limits_{\Gamma_{1}}\vert( \widetilde{v}_{n}+\phi)_{t}\vert^{m-2}
(\widetilde{v}_{ttn}+\phi_{tt}) ^{2}d\sigma = \frac{4}{m^{2}}\displaystyle \int\limits_{\Gamma_{1}}
\left(\frac{\partial}{\partial t}\left(\vert\widetilde{v}_{tn}(t) + \phi_t\vert^{\frac{m-2}{2}}
(\widetilde{v}_{tn}(t)+\phi_t \right) \right)^{2} d\sigma .
$$
Now, integrating equation (\ref{diffsyst}) over $(0,t) $, using the inequality (\ref{estim_0}) and Young's 
and Poincar\'{e}'s inequalities (as in (\ref{lhs4})), there exists $\widetilde{C}_{T} > 0$ such that:
\begin{eqnarray*}
&&\frac{1}{2}\Bigl[ \Vert \nabla \widetilde{v}_{tn}(t) \Vert_{2}^{2}+\Vert \widetilde{v}_{ttn}(t)\Vert_{2}^{2}+
\Vert \widetilde{v}_{ttn}(t)\Vert_{2,\Gamma_{1}}^{2}\Bigl] + 
\alpha \int_{0}^{t}\Vert \nabla \widetilde{v}_{ttn}(s)\Vert_{2}^{2}+ \nonumber \\
&&+ \frac{4r\left( m-1\right) }{m^{2}}\displaystyle \int\limits_{\Gamma_{1}}
\left(\frac{\partial}{\partial t}\left(\vert\widetilde{v}_{tn}(t) + \phi_t\vert^{\frac{m-2}{2}}
(\widetilde{v}_{tn}(t)+\phi_t \right) \right)^{2} d\sigma  \leq \widetilde{C}_{T}.
\end{eqnarray*}
Consequently,we deduce the following results:
\begin{equation}\label{regul3}
\begin{array}{lll}
\left(\widetilde{v}_{ttn}(t)\right)_{n\in \mathbb{N}} &\text{ is bounded in }&
L^{\infty }\left(0,T;L^{2}(\Omega) \right) \quad,\\
\left(\widetilde{v}_{ttn}(t)\right)_{n\in \mathbb{N}} &\text{ is bounded in }&
L^{\infty }\left(0,T;L^{2}\left( \Gamma_{1}\right) \right)  \quad,\\
\left(\widetilde{v}_{tn}(t)\right)_{n\in \mathbb{N}} &\text{ is bounded in }&
L^{\infty }\left(0,T;H_{\Gamma_{0}}^{1}(\Omega) \right) \quad . 
\end{array}
\end{equation}
From (\ref{regul1}), (\ref{regul1bis}), (\ref{regul2}) and (\ref{regul3}), we have 
$\left(\widetilde{v}_{n}\right)_{n\in \mathbb{N}}$ is bounded in $L^{\infty }\left(0,T;H_{\Gamma_{0}}^{1}(\Omega)\right)$.
Then $\left(\widetilde{v}_{n}\right)_{n\in \mathbb{N}}$ is bounded in $L^{2}\left(0,T;H_{\Gamma_{0}}^{1}(\Omega) \right)$.
Since $\left(\widetilde{v}_{tn}\right)_{n\in \mathbb{N}}$ is bounded in $L^{\infty }\left( 0,T;L^{2}(\Omega)\right)$, 
$\left(\widetilde{v}_{tn}\right)_{n\in \mathbb{N}}$ is bounded in $L^{2}\left(0,T;L^{2}(\Omega)\right)$.
Consequently $\left(\widetilde{v}_{n}\right)_{n\in \mathbb{N}}$ is bounded in 
$H^{1}\left( 0,T;H^{1}(\Omega) \right)$. \\
Since the embedding $H^{1}\left(0,T;H^{1}(\Omega)\right) \hookrightarrow L^{2}\left(0,T;L^{2}(\Omega) \right)$ 
is compact, by using Aubin-Lions theorem, we can extract a subsequence $\left(\widetilde{v}_{\mu}\right)_{\mu \in \mathbb{N}}$ 
of $\left(\widetilde{v}_{n}\right)_{n\in \mathbb{N}}$ such that
$$
\widetilde{v}_{\mu }\rightarrow \widetilde{v}\text{ strongly in }L^{2}\left(0,T;L^{2}(\Omega) \right) .
$$
Therefore,
$$
\widetilde{v}_{\mu }\rightarrow \widetilde{v}\ \text{strongly and a.e on } \left( 0,T\right) \times \Omega \quad .
$$
Following \cite[Lemme 3.1]{L69}, we get:
$$
\vert \widetilde{v}_{\mu }\vert^{p-2}\widetilde{v}_{\mu}\rightarrow \vert \widetilde{v}\vert^{p-2}\widetilde{v}\ 
\text{strongly and a.e on }\left( 0,T\right) \times \Omega .
$$
On the other hand, we already have proved in the preceding section that:
\begin{eqnarray*}
&&\left(\widetilde{v}_{tn}\right)_{n\in\mathbb{N}}\text{ is bounded in }L^{\infty }\left( 0,T;L^{2}\left(\Gamma_{1}\right) \right)
\end{eqnarray*}%
From (\ref{regul1}) and (\ref{regul3}), since
$$
\Vert \widetilde{v}_{n}(t) \Vert_{H^{\frac{1}{2}}(\Gamma_1)} \leq C \Vert \nabla \widetilde{v}_{n}(t) \Vert_{2}
\text{ and }
\Vert \widetilde{v}_{tn}(t) \Vert_{H^{\frac{1}{2}}(\Gamma_1)} \leq C \Vert \nabla \widetilde{v}_{tn}(t) \Vert_{2}
$$
we deduce that:
\begin{eqnarray*}
\left(\widetilde{v}_{n}\right)_{n\in\mathbb{N}}&\text{ is bounded in }&L^{2}\left( 0,T;H^{\frac{1}{2}}(\Gamma_1)\right) \\
\left(\widetilde{v}_{tn}\right)_{n\in\mathbb{N}}&\text{ is bounded in }&L^{2}\left( 0,T;H^{\frac{1}{2}}(\Gamma_1)\right) 
\quad .\\
\left(\widetilde{v}_{ttn}\right)_{n\in\mathbb{N}}&\text{ is bounded in }&L^{2}\left( 0,T;L^2(\Gamma_1)\right)
\end{eqnarray*}
Since the embedding $H^{\frac{1}{2}}(\Gamma_{1}) \hookrightarrow L^{2}(\Gamma_{1})$ is compact,
again by using Aubin-Lions theorem, we conclude that we can extract a subsequence also denoted $\left(\widetilde{v}_{\mu}\right)_{\mu \in \mathbb{N}}$ of $\left(\widetilde{v}_{n}\right)_{n\in \mathbb{N}}$ such that:
\begin{equation} \label{strong1}
 \widetilde{v}_{t\mu} \rightarrow \widetilde{v}_t \mbox{ strongly in }L^{2}\left( 0,T;L^2(\Gamma_1)\right) \quad .  
\end{equation}
Therefore from (\ref{regul2}), we obtain that:
$$
\vert \widetilde{v}_{t\mu }\vert^{m-2}\widetilde{v}_{t\mu} \rightharpoonup \chi \text{ weakly in }
L^{m'}\left(\left(0,T\right) \times \Gamma_{1}\right) ,
$$
It suffices to prove now that $\chi =\vert \widetilde{v}_{t}\vert^{m-2}\widetilde{v}_{t}$.\\
Clearly, from (\ref{strong1}) we get:
$$
\vert \widetilde{v}_{t\mu }\vert^{m-2}\widetilde{v}_{t\mu}\rightarrow \vert \widetilde{v}_t\vert^{m-2}\widetilde{v}_t\ 
\text{strongly and a.e on }\left( 0,T\right) \times \Gamma_{1}.
$$
Again, by using the Lions's lemma, \cite[Lemme 1.3]{L69}, we obtain $\chi =\vert \widetilde{v}_{t}\vert^{m-2}\widetilde{v}_{t}.$
The proof now can be completed arguing as in \cite[Th\'eor\`eme 3.1]{L69}.
\subsubsection{Uniqueness}
Let $v,w$ two solutions of the problem (\ref{ondes_f}) which share the same initial data. Let us denote $z=v-w.$ 
It is staightforward to see that $z$ satisfies:
\begin{equation}\label{eqforz}
\left.
\begin{array}{c}
\Vert \nabla z\Vert_{2}^{2}+\Vert \nabla z_{t}\Vert_{2}^{2}+\Vert z_{t}\Vert_{2,\Gamma_{1}}^{2}+
2\alpha\displaystyle \int\limits_{0}^{t}\Vert \nabla z\Vert_{2}^{2}ds \\
+2r\displaystyle \int\limits_{0}^{t}\displaystyle \int\limits_{\Gamma_{1}}
\Bigl[ \vert v_{t}(s) \vert^{m-2}v_{t}(s) -\vert w_{t}(s) \vert^{m-2}w_{t}(s) \Bigl] (v_t(s) - w_t(s) ) ds d\sigma = 0 \quad .
\end{array}%
\right.
\end{equation}%
By using the algebraic inequality: 
\begin{equation}\label{algebr_m}
\forall \, m \geq 2 \,,\, \exists \, c > 0 \,,\, \forall \, a\,,\, b\in \mathbb{R} \;,\;
\Bigl[ \vert a\vert^{m-2}a-\vert b\vert^{m-2}b\Bigl](a-b) \, \geq \, c \, \vert a-b\vert^{m}
\end{equation}
equation (\ref{eqforz}) yields to:
$$
\left.
\begin{array}{c}
\Vert z_{t}\Vert_{2}^{2}+\Vert \nabla z\Vert_{2}^{2}+\Vert z_{t}\Vert_{2,\Gamma_{1}}^{2}+
2\alpha\displaystyle \int\limits_{0}^{t}\Vert \nabla z_{t}\Vert_{2}^{2}ds \\
+c\displaystyle \int\limits_{0}^{t}\displaystyle \int\limits_{\Gamma_{1}}\vert v_{t}\left(s\right) -w_{t}(s) \vert^{m}
ds d\sigma \leq 0 \quad .
\end{array}
\right.
$$
Then, this last inequality yields to  $z=0$.

This finishes the proof of the lemma \ref{existence_f}.
\subsection{Proof of lemma \ref{existence_u}}

We first approximate $u\in C\left([ 0,T] ,H_{\Gamma_{0}}^{1}(\Omega)\right) \cap C^{1}\left([ 0,T] ,L^{2}(\Omega)\right)$ 
endowed with the standard norm 
$\Vert u \Vert = \displaystyle \max_{t \in [0,T]}  \Vert u_t (t)\Vert_2 + \Vert u(t )\Vert_{H^1(\Omega)}$, by  a sequence  $(u^{k})_{k\in\mathbb{N}} \subset C^{\infty}([ 0,T] \times \overline{\Omega})$ 
by a standard convolution arguments (see \cite{B83}). It is clear that $f\left( u^{k}\right) =\vert u^{k}\vert^{p-2}u^{k}\in H^{1}\left( 0,T;L^{2}(\Omega) \right)$.
This type of approximation has been already used by Vitillaro in \cite{V02_1, V02_2}.
Next, we approximate the initial data $u_{1} \in L^{2}(\Omega)$ by a sequence $(u_{1}^{k})_{k\in\mathbb{N}} \subset C_{0}^{\infty}(\Omega)$.
 Finally, since the space $H^2(\Omega) \cap V \cap H_{\Gamma_{0}}^{1}(\Omega)$ is dense in $ H_{\Gamma_{0}}^{1}(\Omega)$ for the $H^1$ norm, we approximate
$u_0 \in H_{\Gamma_{0}}^{1}(\Omega)$ by a sequence $(u_{0}^{k})_{k\in\mathbb{N}} \subset H^2(\Omega) \cap V \cap H_{\Gamma_{0}}^{1}(\Omega)$.

We consider now the set of the following problems:
\begin{equation}\label{approx_k}
\left\{
\begin{array}{ll}
\hspace*{-0.3cm} v_{tt}^k-\Delta v^k-\alpha \Delta v_{t}^k=\vert u^k\vert^{p-2}u^k, & x\in \Omega ,\ t>0 \\[0.1cm]
\hspace*{-0.3cm}v^k(x,t) =0, &  x\in \Gamma_{0},\ t>0  \\[0.1cm]
\hspace*{-0.3cm}v_{tt}^k(x,t) =- \left[\displaystyle \frac{\partial v^k}{\partial \nu }(x,t) +
\frac{\alpha \partial v_{t}^k}{\partial \nu }(x,t) +r\vert v_{t}^k\vert^{m-2}v_{t}^k( x,t) \right]
&  x\in \Gamma_{1},\ t>0 \\[0.1cm]
\hspace*{-0.3cm}v^k(x,0) = u^k_0, \; v_{t}^k( x,0)= u^k_1  & x\in \Omega \quad .
\end{array}
\right.
\end{equation}
Since every hypothesis of  lemma \ref{existence_f} are verified, we can find a sequence of unique solution $\left( v_{k}\right)_{k \in \mathbb{N}}$
of the problem (\ref{approx_k}). Our goal now is to show that $(v^{k},v_t^k)_{k\in \mathbb{N}}$ is a Cauchy sequence 
in the space
$$ \begin{array}{lll}
Y_{T} &= \Bigl\{ &(v,v_{t}) /v\in C\left( \left[ 0,T\right],H_{\Gamma_{0}}^{1}(\Omega) \right) 
\cap C^{1}\left( \left[ 0,T\right] ,L^{2}(\Omega) \right) ,  \\
&&v_{t}\in L^{2}\left( 0,T;H_{\Gamma_{0}}^{1}(\Omega)
\right) \cap L^{m}\left( \left( 0,T\right) \times \Gamma_{1}\right)\Bigl\}
\end{array}
$$
endowed with the norm
$$
\Vert (v,v_t)\Vert_{Y_{T}}^2=\max_{0 \leq t \leq T}\Bigl[ \Vert v_{t}\Vert_{2}^2 + 
\Vert \nabla v \Vert_{2}^2\Bigl] + \Vert v_{t} \Vert^2_{L^{m}\bigl( \left( 0,T\right) \times \Gamma_{1}\bigl) } +
\int_0^t \Vert \nabla v_t(s) \Vert_2^2 \; ds \;  .
$$
For this purpose, we set $U=u^{k}-u^{k'},\ V=v^{k}-v^{k'}$ .
It is straightforward to see that $V$ satisfies:
\begin{equation*}
\left\{
\begin{array}{lll}
\hspace*{-0.2cm} V_{tt}-\Delta V & -\alpha \Delta V_{t}=\vert u^k\vert^{p-2} u^k - 
\vert u^{k'} \vert^{p-2} u^{k'}  &x\in \Omega ,\ t>0 \\[0.1cm]
\hspace*{-0.2cm}V(x,t)  =& 0 &  x\in \Gamma_{0},\ t>0  \\[0.1cm]
\hspace*{-0.2cm}V_{tt}(x,t) =&- \Bigl[\displaystyle \frac{\partial V}{\partial \nu }(x,t) +
\frac{\alpha \partial V_{t}}{\partial \nu }(x,t) \Bigl] - & \\[0.3cm]
&r \Bigl(\vert v_{t}^k\vert^{m-2}v_{t}^k( x,t) - \vert v_{t}^{k'} \vert^{m-2}v_{t}^{k'}( x,t) \Bigl)
&  x\in \Gamma_{1},\ t>0 \\[0.1cm]
\hspace*{-0.2cm}V(x,0) = &u_0^k - u_0^{k'}, \; V_{t}( x,0)= u_1^k - u_1^{k'}  & x\in \Omega \quad .
\end{array}
\right.
\end{equation*}
We multiply the above differential equations by $V_{t}$, we integrate over $(0,t) \times \Omega$ and we use
integration by parts to obtain:
$$
\left.
\begin{array}{ll}
\displaystyle\frac{1}{2}&\left( \Vert V_{t}\Vert_{2}^{2}+\Vert \nabla V\Vert_{2}^{2}+\Vert V_{t}\Vert_{2,\Gamma_{1}}^{2}\right)
+\alpha \displaystyle \int\limits_{0}^{t}\Vert \nabla V_{t}\Vert_{2}^{2}\, ds \\
+&r\displaystyle \int\limits_{0}^{t}\displaystyle \int\limits_{\Gamma_{1}}
\left( \vert v_{t}^{k}\vert^{m-2} v_{t}^{k} - \vert v_{t}^{k'}\vert^{m-2} v_{t}^{k'}\right) 
\left(v_{t}^{k}-v_{t}^{k'}\right) d\sigma ds \\
=&\displaystyle \frac{1}{2}
\left( \Vert V_{t}(0) \Vert _{2}^{2}+\Vert \nabla V(0) \Vert_{2}^{2}+\Vert V_{t}(0)\Vert_{2,\Gamma_{1}}^{2}\right)  \\
+&\displaystyle \int\limits_{0}^{t}\displaystyle \int\limits_{\Omega}
\left( \vert u^{k}\vert^{p-2}u^{k}-\vert u^{k'}\vert^{p-2}u^{k'}\right) \left( v_{t}^{k}-v_{t}^{k'}\right) \, dx d\tau 
,\quad \forall t\in \left( 0,T\right) \quad .
\end{array}
\right. 
$$
By using the algebraic inequality (\ref{algebr_m}), we get:
\begin{eqnarray*}
&\displaystyle \frac{1}{2}
& \left(\Vert V_{t}\Vert_{2}^{2}+\Vert \nabla V\Vert_{2}^{2}+\Vert V_{t}\Vert_{2,\Gamma_{1}}^{2}\right) +
\alpha \displaystyle \int\limits_{0}^{t}\Vert \nabla V_{t}\Vert_{2}^{2}ds+c_{1}\Vert V_{t}\Vert _{m,\Gamma_{1}}^{m} \\
&\leq &\frac{1}{2}
\left( \Vert V_{t}(0) \Vert_{2}^{2}+\Vert \nabla V(0) \Vert_{2}^{2}+\Vert V_{t}(0) \Vert_{2,\Gamma_{1}}^{2}\right) \\
&& + \displaystyle \int\limits_{0}^{t}\displaystyle \int\limits_{\Omega}
\left( \vert u^{k}\vert^{p-2}u^{k}-\vert u^{k'}\vert^{p-2}u^{k'}\right) \left( v_{t}^{k}-v_{t}^{k'}\right) \, dxd\tau ,\quad
 \forall t\in \left( 0,T\right) .
\end{eqnarray*}
In order to find a majoration of the term:
$$
\displaystyle \int\limits_{0}^{t}\displaystyle \int\limits_{\Omega}
\left( \vert u^{k}\vert^{p-2}u^{k}-\vert u^{k'}\vert^{p-2}u^{k'}\right) \left( v_{t}^{k}-v_{t}^{k'}\right) \, dxd\tau ,\quad
 \forall t\in \left( 0,T\right)
$$
in the previous inequality, we use the result of Georgiev and Todorova \cite{GT94} 
(specifically their equations (2.5) and (2.6) in proposition 2.1).
The hypothesis on $p$ ensures us to use exactly the same argument.
Thus by applying Young's inequality and Gronwall inequality, there exists $C$ depending only on $\Omega \mbox{ and }p$ 
such that:
$$
\Vert V\Vert_{Y_{T}}\leq 
C\left( \Vert V_{t}(0) \Vert_{2}^{2}+\Vert \nabla V(0) \Vert _{2}^{2}+\Vert V_{t}(0) \Vert_{2,\Gamma_{1}}^{2}\right) +
CT\Vert U\Vert _{Y_{T}}.
$$
Let us now remark that from the notations used above, we have:
$$ V(0) =  u_0^{k}-u^{k'}_0 \mbox{ , } V_{t}(0) = u_1^{k}-u^{k'}_1 \mbox{ and } U =  u^{k}-u^{k'} \quad .$$
Thus, since $(u_{0}^{k})_{k\in\mathbb{N}} $ is a converging sequence in $H_{\Gamma_{0}}^{1}\left(\Omega \right)$, 
$( u_{1}^{k})_{k\in\mathbb{N}} $ is a converging sequence in $L^{2}\left(\Omega \right) $ and 
$\left( u^{k}\right)_{k\in\mathbb{N}} $ is a converging sequence in 
$C\left( \left[ 0,T\right] ,H_{\Gamma_{0}}^{1}(\Omega) \right) \cap 
C^{1}\left(\left[ 0,T\right] ,L^{2}(\Omega) \right) $ (so in $Y_{T}$ also), we conclude that 
$(v^{k},v_t^k)_{k\in\mathbb{N}} $ is a Cauchy sequence in $Y_{T}$. 
Thus $(v^{k},v_t^k) $ converges to a limit $(v,v_t) \in Y_{T}$. Now  by the same procedure used by Georgiev and
Todorova in \cite{GT94}, we prove that this limit is a weak solution of the problem (\ref{ondes_u}).
This completes the proof of the lemma \ref{existence_u}.
\subsection{Proof of theorem \ref{existence}}
In order to prove theorem \ref{existence}, we use the contraction mapping theorem.\\ 
For $T>0,$ let us define the convex closed subset of $Y_T$:
$$
X_{T}= \left\{ (v,v_t) \in Y_T \mbox{ such that } v(0) =u_{0},v_{t}(0) =u_{1}\right\} \quad .
$$
Let us denote:
$$
B_{R}\left( X_{T}\right) =\left\{ v\in X_{T};\Vert v\Vert_{Y_{T}}\leq R\right\} \quad .
$$
Then, lemma \ref{existence_u} implies that for any $u\in X_{T}$, we may define 
$v=\Phi\left( u\right) $ the unique solution of (\ref{ondes_u}) corresponding to $u$.
Our goal now is to show that for a suitable $T>0$, $\Phi $ is a contractive map satisfying 
$\Phi \left( B_R(X_{T})\right) \subset B_R(X_{T})$ . \\
Let $u\in B_R(X_{T})$ and $v=\Phi \left( u\right) $. Then for all $t \in [ 0,T] $ we have:
\begin{equation}\label{schauder}
\begin{array}{l}
\Vert v_{t}\Vert_{2}^{2}+\Vert \nabla v\Vert_{2}^{2}+\Vert v_{t}\Vert_{2,\Gamma_{1}}^{2}+
2\displaystyle \int\limits_{0}^{t}\Vert v_{t}\Vert _{m,\Gamma_{1}}^{m}ds +
2\alpha\displaystyle \int\limits_{0}^{t}\Vert \nabla v_{t}\Vert_{2}^{2} \, ds\\
=\Vert u_{1}\Vert_{2}^{2}+\Vert \nabla u_{0}\Vert_{2}^{2}+\Vert u_{1}\Vert_{2,\Gamma_{1}}^{2}+
2\displaystyle \int\limits_{0}^{t}\displaystyle \int\limits_{\Omega}
\vert u\left( \tau\right) \vert^{p-2}u\left( \tau \right) v_{t}\left( \tau \right) \, dx \, d\tau
\end{array} \quad .
\end{equation}
Using H\"older inequality, we can control the last term in the right hand side of the inequality (\ref{schauder})
as follows:
$$
\displaystyle \int\limits_{0}^{t}\displaystyle \int\limits_{\Omega}\vert u\left( \tau \right)
\vert^{p-2}u\left( \tau \right) v_{t}\left( \tau \right) dxd\tau
\leq  \displaystyle \int\limits_{0}^{t}\Vert u\left( \tau \right) \Vert_{2N/\left( N-2\right) }^{p-1}
\Vert v_{t}\left( \tau \right)\Vert_{{2N}/{\bigl(3N - N p + 2 (p-1)\bigl)} } d\tau
$$
Since $ \displaystyle p \leq \frac{2 N }{N-2}$, we have 
$
\displaystyle \frac {2N}{\bigl(3N - N p + 2 (p-1)\bigl)} \leq \frac{2 N } {N-2} \quad .
$\\[0.1cm]
Thus, by Young's and Sobolev's inequalities, we get $\forall \,\delta > 0 \,,\, \exists \, C(\delta) > 0$, such that:
$$
\forall t\in \left(0,T\right)\,,\, 
\; \displaystyle \int\limits_{0}^{t}\displaystyle \int\limits_{\Omega}\vert u\left( \tau \right)
\vert^{p-2}u\left( \tau \right) v_{t}\left( \tau \right) dxd\tau
\leq  C(\delta) t R^{2(p-1)}+ \delta \displaystyle \int\limits_{0}^{t}
\Vert \nabla v_{t}\left( \tau \right) \Vert_{2}^{2}d\tau .
$$
Inserting the last estimate in the inequality  (\ref{schauder})  and choosing $\delta$ small enough in order to 
counter-balance the last term of the left hand side of the inequality (\ref{schauder}) we get:
$$
\Vert v\Vert _{Y_{T}}^{2}\leq \frac{1} {2} R^{2} +C T R^{2(p-1) }.
$$
Thus, for $T$ sufficiently small, we have $\Vert v\Vert_{Y_{T}}\leq R$. This shows that $v\in B_{R}\left( X_{T}\right)$.\\
Next, we have to verify that $\Phi $ is a contraction.
To this end, we set $U~=~u~-~\bar{u}$ and $V~=~v~-~\bar{v}$, 
where $ v = \Phi(u) $ and $\bar{v}=\Phi(\bar{u})$ are the solutions of problem (\ref{ondes_u}) corresponding respectively
to $u $ and $v$. Consequently we have: 
\begin{equation}\label{eqcontract}
\left\{
\begin{array}{lll}
\hspace*{-0.2cm} V_{tt}-\Delta V & -\alpha \Delta V_{t}=\vert u\vert^{p-2} u - 
\vert \bar{u} \vert^{p-2} \bar{u}  &x\in \Omega ,\ t>0 \\[0.1cm]
\hspace*{-0.2cm}V(x,t)  =& 0 &  x\in \Gamma_{0},\ t>0  \\[0.1cm]
\hspace*{-0.2cm}V_{tt}(x,t) =&- \Bigl[\displaystyle \frac{\partial V}{\partial \nu }(x,t) +
\frac{\alpha \partial V_{t}}{\partial \nu }(x,t) \Bigl] - & \\[0.3cm]
&r \Bigl(\vert v_{t}\vert^{m-2}v_{t}( x,t) - \vert \bar{v}_{t} \vert^{m-2}\bar{v}_{t}(x,t) \Bigl)
&  x\in \Gamma_{1},\ t>0 \\[0.1cm]
\hspace*{-0.2cm}V(x,0) = &0, \; V_{t}( x,0)= 0& x\in \Omega \quad .
\end{array}
\right.
\end{equation}
By multiplying the differential equation (\ref{eqcontract}) by $V_{t}$ and integrating over 
$(0,t) \times \Omega$, we get:
$$
\begin{array}{ll}
\displaystyle \frac{1}{2}&\left(\Vert V_{t}\Vert_{2}^{2}+\Vert \nabla V\Vert_{2}^{2}+\Vert V_{t}\Vert_{2,\Gamma_{1}}^{2}\right)
+\alpha \displaystyle \int\limits_{0}^{t}
\Vert \nabla V_{t}\Vert_{2}^{2}ds+\\
& r\displaystyle \int\limits_{0}^{t}\displaystyle \int\limits_{\Gamma_{1}}
\left( \vert v_{t}\vert^{m-2}v_{t} - \vert \bar{v}_{t}\vert^{m-2}\bar{v}_{t}\right)(v_{t}-\bar{v}_{t}) d\sigma ds \\
&=\displaystyle \int\limits_{0}^{t}\displaystyle \int\limits_{\Omega}
\left( \vert u\vert ^{p-2}u-\vert \bar{u}\vert^{p-2}\bar{u}\right) (v_{t}-\bar{v}_{t}) dx d\tau ,\ 
\forall t\in( 0,T)
\end{array} \quad .
$$
Again, by using the algebraic inequality (\ref{algebr_m}),  we have:
\begin{equation}\label{contract1}
\begin{array}{l}
\displaystyle \frac{1}{2}\left(\Vert V_{t}\Vert_{2}^{2}+\Vert \nabla V\Vert_{2}^{2}+\Vert V_{t}\Vert_{2,\Gamma_{1}}^{2}\right) +
\alpha \displaystyle \int\limits_{0}^{t}\Vert \nabla V_{t}\Vert_{2}^{2}ds+c_{1}\Vert V_{t}\Vert _{m,\Gamma_{1}}^{m} \\
\leq \displaystyle \int\limits_{0}^{t}\displaystyle \int\limits_{\Omega}
\left( \vert u\vert^{p-2}u-\vert \bar{u}\vert^{p-2}\bar{u}\right)(v_{t}-\bar{v}_{t}) dx d\tau 
,\ \forall t\in \left( 0,T\right)
\end{array}
\quad .
\end{equation}
To estimate the term in the right hand side of the inequality (\ref{contract1}), let us denote:
$$
I(t) :=\displaystyle \int\limits_{0}^{t}\displaystyle \int\limits_{\Omega}
\left(\vert u\vert^{p-2}u-\vert \bar{u}\vert^{p-2}\bar{u}\right)(v_{t}-\bar{v}_{t}) dxd\tau \quad .
$$
Using the algebraic inequality:
$$
\left \vert \vert u\vert^{p-2}u-\vert \bar{u}\vert^{p-2}\bar{u}\right\vert \leq 
c_{p}\vert u-\bar{u}\vert\left( \vert u\vert^{p-2}+\vert \bar{u}\vert^{p-2}\right) ,
$$
which holds for any $u,\bar{u}\in \mathbb{R}$, where $c_{p}$ is a positive constant depending only on $p$, we find:
$$
I(t) \leq c_{p}\displaystyle \int\limits_{0}^{T}\displaystyle \int\limits_{\Omega}
\vert u-\bar{u}\vert\left( \vert u\vert^{p-2}+\vert \bar{u}\vert^{p-2}\right) \vert V_{t}\vert dx d\tau \quad .
$$
Following the same argument as Vitillaro in \cite[eq 77]{V02_2}, choosing $p<r_{0}<\bar{q}$ such that: 
$$\frac{\bar{q}}{\bar{q}-p+1}<\frac{r_{0}}{r_{0}-p+1}<m \quad ,$$ 
let $s>1$ such that: 
$$
\frac{1}{m}+\frac{1}{r_{0}}+\frac{1}{s}=1 \quad .
$$ 
Using H\"older's inequality we obtain:
\begin{equation}\label{eqint1}
I(t) \leq c_{p}\displaystyle \int\limits_{0}^{T}\left( \Vert u-\bar{
u}\Vert _{r_{0}}\Vert V_{t}\Vert _{m}\right) .\left(
\displaystyle \int\limits_{\Omega}\left( \vert u\vert^{p-2}+\vert
\bar{u}\vert^{p-2}\right) ^{s}\right) ^{1/s}.
\end{equation}
Therefore, the algebraic inequality  (\ref{algebr_lambda}) gives us:
$$
\left( \displaystyle \int\limits_{\Omega}
\left(\vert u \vert^{p-2} + \vert \bar{u} \vert^{p-2} \right)^{s}\right)^{1/s}\leq 
2^{s-1} \left( \Vert u \Vert_{(p-2)s}^{(p-2) s} + \Vert \bar{u}\Vert_{(p-2) s}^{( p-2) s}\right) ^{1/s} \quad .
$$
But since
$$
\left( A+B\right) ^{\beta }\leq A^{\beta }+B^{\beta }\,,\, \forall \, A\,,\,B\geq 0 \text{ and } 0< \beta <1
$$
we get
\begin{equation} \label{eqint2}
\left( \displaystyle \int\limits_{\Omega}
\left( \vert u\vert^{p-2}+\vert \bar{u}\vert^{p-2}\right)^{s}\right)^{1/s} \leq 
2^{s-1}\left( \Vert u\Vert_{(p-2)s}^{(p-2)} + \Vert \bar{u}\Vert_{(p-2) s}^{(p-2)}\right) 
\quad .
\end{equation}
Consequently, inserting the inequality (\ref{eqint1}) in (\ref{eqint2}) and using
Poincar\'e's inequality, we obtain:
$$
I(t) \leq c_{2}R^{p-2}\displaystyle \int\limits_{0}^{T}\Vert u-\bar{u}\Vert _{r_{0}}\Vert \nabla V_{t}\Vert_{2}ds.
$$
Applying H\"{o}lder's inequality in time, we finally get:
\begin{equation}\label{eqint3}
\begin{array}{ll}
I(t)  &\leq c_{2}R^{p-2}T^{1/2}\Vert u-\bar{u}\Vert_{L^{\infty }\left( 0,T;L^{r_{0}}(\Omega) \right)} \; 
\left(\displaystyle \int\limits_{0}^{T}\Vert \nabla V_{t}\Vert_{2}^{2}\right) ^{1/2}  \\
&\leq \displaystyle \frac{c_{2}}{2}R^{p-2}T^{1/2}
\left[ 
\Vert u-\bar{u}\Vert _{L^{\infty }(0,T;L^{r_{0}}(\Omega))}^{2} +
\displaystyle \int\limits_{0}^{T}\Vert \nabla V_{t}\Vert_{2}^{2}
\right] \; .
\end{array}
\end{equation}
Lastly, by choosing $T$  small enough in order to have: 
$$\alpha -\frac{c_{2}}{2}R^{p-2}T^{1/2}>0 \quad ,$$ 
we conclude by inserting the estimate (\ref{eqint3}) in the  estimate (\ref{contract1}) that:
\begin{equation}\label{contract2}
\begin{array}{ll}
&\displaystyle \frac{1}{2}\left(\Vert V_{t}\Vert_{2}^{2}+\Vert \nabla V\Vert_{2}^{2}+\Vert V_{t}\Vert_{2,\Gamma_{1}}^{2}\right)
+\alpha \displaystyle \int\limits_{0}^{t}\Vert \nabla V_{t}\Vert_{2}^{2}ds+c_{1}\Vert V_{t}\Vert _{m,\Gamma_{1}}^{m}\\
&\hspace*{2cm} 
\leq \displaystyle \frac{c_{2}}{2}R^{p-2}T^{1/2}\Vert u-\bar{u}\Vert_{L^{\infty }\left( 0,T;L^{r_{0}}(\Omega) \right) }^{2} 
\quad .
\end{array}
\end{equation}
Since $r_{0}<\bar{q},$ using the embedding
$$
L^{\infty }\left(0,T;H_{\Gamma_{0}}^{1}(\Omega) \right) \hookrightarrow L^{\infty }\left( 0,T;L^{r_{0}}(\Omega) \right)
$$
in the estimate (\ref{contract2}), we finally have:
\begin{equation}\label{contract3}
\Vert V\Vert_{Y_{T}}^{2}\leq c_{3}R^{p-2}T^{1/2} \Vert U\Vert_{Y_{T}}^{2}.
\end{equation}%
By choosing $T$ small enough in order to have  
$$c_{3 }R^{p-2} T^{1/2} < 1 \quad .$$ 
the estimate (\ref{contract3}) shows that $\Phi $ is a contraction. 
Consequently the contraction mapping theorem guarantees the existence of a unique $v$ satisfying $v=\Phi(v)$. 
The proof of theorem \ref{existence} is now completed.
\begin{remark}\rm
To prove the existence and uniqueness of the solution to the more general problem:
\begin{equation*}
\left\{
\begin{array}{ll}
u_{tt}-\Delta u-\alpha \Delta u_{t}=f(u), & x\in \Omega ,\ t>0 \\[0.1cm]
u( x,t) =0, &  x\in \Gamma_{0},\ t>0  \\[0.1cm]
u_{tt}(x,t) =- \left[\displaystyle \frac{\partial u}{\partial \nu }(x,t) +
\frac{\alpha \partial u_{t}}{\partial \nu }(x,t) +g(u_t) \right]
&  x\in \Gamma_{1},\ t>0 \\[0.1cm]
u( x,0) =u_{0}(x), \; u_{t}( x,0) =u_{1}(x) & x\in \Omega \quad .
\end{array}
\right.
\end{equation*}
we can use the same method, provided that the functions $f$ and $g$ satisfy respectively the  conditions
$(H_3) - (H_7)$ and $(H_8) - (H_9)$ of the paper of Calvacanti et \textit{al.} \cite{CCSM00}.
\end{remark}
\section{Exponential growth}
In this section we consider the  problem (\ref{ondes}) and we will prove that when the initial
data are large enough (in the energy point of view), the energy grows exponentially and thus so the $L^p$ norm.\\
In order to state and prove the result, we introduce the following notations.
Let $B$ be the best constant of the embedding $H_{0}^{1}(\Omega) \hookrightarrow L^{p}(\Omega) $ defined by:
$$
B^{-1} = \mbox{inf}\left\{\Vert \nabla u \Vert_2 : u \in  H_{0}^{1}(\Omega), \Vert u\Vert_p = 1 \right\} \quad .
$$
We also define the energy functional:
\begin{equation}\label{energy}
E(u(t)) = E(t) =\frac{1}{2}\left \Vert \nabla u \right \Vert _{2}^{2} -
\frac{1}{p}\left\Vert u\right\Vert _{p}^{p} +
\frac{1}{2}\left\Vert u_{t}\right\Vert _{2}^{2}+\frac{1}{2} \left\Vert u_{t}\right\Vert
_{2,\Gamma _{1}}^{2}.
\end{equation}
Finally we define the following constant which will play an important role in the proof of our result:
\begin{equation}\label{constant}
\alpha_{1} = B^{-p/( p-2)} \,,\mbox{ and } d=(\frac{1}{2}-\frac{1}{p}) \alpha _{1}^{2}.
\end{equation}
In order to obtain the exponential growth of the energy, we will use the following lemma (see Vitillaro \cite{V99}, for the proof):

\begin{lemma}\label{Vitillaro}
Let  $u$  be a classical solution of (\ref{ondes}). Assume that  
$$ E(0) < d \mbox{ and } \left\Vert \nabla u_{0}\right\Vert _{2} >\alpha _{1}. $$ 
Then there exists a constant $\alpha _{2}>\alpha _{1}$ such that 
\begin{equation} \label{estimateH1}
\left\Vert \nabla u( .,t) \right\Vert_2 \geq \alpha _{2},\ \
\forall t\geq 0,
\end{equation}
and
\begin{equation}\label{estimateLp}
\left\Vert u\right\Vert _{p}\geq B\alpha _{2},\ \ \forall t\geq 0.
\end{equation}%
\end{lemma}
Let us now state our new result.
\begin{theorem} Assume that $m < p$  where  $2<p\leq \bar{q}$. Suppose that 
$$
E(0) <d \mbox{ and } \left\Vert \nabla u_{0}\right\Vert _{2}>\alpha _{1} .
$$
Then the solution of problem (\ref{ondes}) growths exponentially in the $L^{p}$ norm.

\end{theorem}
\textbf{Proof:} By setting
\begin{equation}
H(t) = d -E(t)
\end{equation}%
we get from the definition of the energy (\ref{energy}):
\begin{equation}\label{ineqH}
0< H(0) \leq H(t) \leq d -\left[ \frac{1}{2}\left\Vert u_{t}\right\Vert_{2}^{2}
+\frac{1}{2}\left\Vert u_{t}\right\Vert _{2,\Gamma_{1}}^{2}
+\frac{1}{2}\left\Vert \nabla u\right\Vert _{2}^{2}-\frac{1}{p}
\left\Vert u\right\Vert _{p}^{p}\right] ,
\end{equation}
using the fundamental estimate (\ref{estimateH1}) and the equality (\ref{constant}), we get:
\[
d-\frac{1}{2}\left\Vert \nabla u\right\Vert _{2}^{2}<d-\frac{1}{2}\alpha
_{1}^{2}=-\frac{1}{p}\alpha _{1}^{2}<0,~\forall t\geq 0.
\]
Hence we finally obtain the following inequality:
\begin{equation*}
0< H( 0) \leq H(t) \leq \frac{1}{p}\left\Vert u \right\Vert _{p}^{p},\quad \forall t \geq 0.
\end{equation*}
For $\varepsilon $ small to be chosen later, we then define the auxillary function:
\begin{equation}\label{defL}
L(t) = H(t) +\varepsilon \int_{\Omega}u_{t}u dx+
\varepsilon \int_{\Gamma _{1}}u_{t}ud\sigma 
+\frac{\varepsilon \alpha }{2}\left\Vert \nabla u\right\Vert _{2}^{2}.
\end{equation}%
Let us remark that $L$ is a small perturbation of the energy. By taking the time derivative of (\ref{defL}),
we obtain:
\begin{eqnarray}\label{derivL}
\frac{dL(t) }{dt} & = &\alpha \left\Vert \nabla u_{t}\right\Vert
_{2}^{2}+r\left\Vert u_{t}\right\Vert _{m,\Gamma _{1}}^{m}+\varepsilon
\left\Vert u_{t}\right\Vert _{2}^{2}  + \varepsilon \alpha \int_{\Omega }\nabla u_{t}\nabla u dx\nonumber \\
&&+\varepsilon \int_{\Omega }u_{tt}udx+\varepsilon \int_{\Gamma
_{1}}u_{tt}ud\sigma +\varepsilon \left\Vert u_{t}\right\Vert _{2,\Gamma
_{1}}^{2}.
\end{eqnarray}
Using problem (\ref{ondes}), the equation (\ref{derivL}) takes the form:
\begin{eqnarray}\label{derivL2}
\frac{dL(t) }{dt} &=&\alpha \left\Vert \nabla u_{t}\right\Vert
_{2}^{2}+r\left\Vert u_{t}\right\Vert _{m,\Gamma _{1}}^{m}+\varepsilon
\left\Vert u_{t}\right\Vert _{2}^{2}-\varepsilon \left\Vert \nabla
u\right\Vert _{2}^{2}  \nonumber \\
&&+\varepsilon \left\Vert u\right\Vert _{p}^{p}+\varepsilon \left\Vert
u_{t}\right\Vert _{2,\Gamma _{1}}^{2}-\varepsilon r\int_{\Gamma
_{1}}\left\vert u_{t}\right\vert ^{m}u_{t}u( x,t) d\sigma.
\end{eqnarray}
To estimate the last term in the right hand side of the previous equality, let $\delta > 0$ be chosen later.
Young's inequality leads to:
\[
\int_{\Gamma _{1}}\left\vert u_{t}\right\vert ^{m}u_{t}u(x,t) d\sigma \leq 
\frac{\delta ^{m}}{m}\left\Vert u\right\Vert _{m,\Gamma _{1}}^{m}+%
\frac{m-1}{m}\delta ^{-m/(m-1) }\left\Vert u_{t}\right\Vert
_{m,\Gamma _{1}}^{m}.
\]%
This yields by substitution in (\ref{derivL2}):
\begin{eqnarray}\label{derivL3}
\frac{dL(t) }{dt} &\geq &\alpha \left\Vert \nabla
u_{t}\right\Vert _{2}^{2}+r\left\Vert u_{t}\right\Vert _{m,\Gamma
_{1}}^{m}+\varepsilon \left\Vert u_{t}\right\Vert _{2}^{2}-\varepsilon
\left\Vert \nabla u\right\Vert _{2}^{2}  \nonumber \\
&&+\varepsilon \left\Vert u\right\Vert _{p}^{p}+\varepsilon \left\Vert
u_{t}\right\Vert _{2,\Gamma _{1}}^{2}-\frac{\varepsilon r}{m}\delta
^{m}\left\Vert u\right\Vert _{m,\Gamma _{1}}^{m}  \\
&&-\frac{\varepsilon r( m-1) }{m}\delta ^{-m/( m-1)
}\left\Vert u_{t}\right\Vert _{m,\Gamma _{1}}^{m} \nonumber  \quad .
\end{eqnarray}%
Let us recall the inequality concerning the continuity of the trace operator
(here and in the sequel, $C$ denotes generic positive constant which may change from line to line):
\[
\left\Vert u\right\Vert _{m,\Gamma _{1}}\leq C\left\Vert u\right\Vert
_{H^{s}(\Omega) },
\]%
which holds for:
$$m\geq 1 \mbox{ and } 0<s<1, s\geq \frac{N}{2}-\frac{N-1}{m}>0$$
and the interpolation and Poincar\'{e}'s inequalities (see \cite{LM68})
\begin{eqnarray*}
\left\Vert u\right\Vert _{H^{s}(\Omega) } &\leq &C\left\Vert
u\right\Vert _{2}^{1-s}\left\Vert \nabla u\right\Vert _{2}^{s} \\
&\leq &C\left\Vert u\right\Vert _{p}^{1-s}\left\Vert \nabla u\right\Vert
_{2}^{s}
\end{eqnarray*}
Thus, we have the following inequality:
\[
\left\Vert u\right\Vert _{m,\Gamma _{1}}\leq C\left\Vert u\right\Vert
_{p}^{1-s}\left\Vert \nabla u\right\Vert _{2}^{s}.
\]%
If $s< 2/m$, using again Young's inequality, we get:
\begin{equation}\label{estiGamma1}
\left\Vert u\right\Vert _{m,\Gamma _{1}}^{m}\leq C\left[ \left( \left\Vert
u\right\Vert _{p}^{p}\right) ^{\frac{m\left( 1-s\right) \mu }{p}}+\left(
\left\Vert \nabla u\right\Vert _{2}^{2}\right) ^{\frac{ms\theta }{2}}\right]
\end{equation}
for $1/\mu +1/\theta =1.$ Here we choose $\theta =2/ms,$ 
to get $\mu =2/\left(2-m s\right)$. Therefore the previous inequality becomes:
\begin{equation}\label{normGamma1}
\left\Vert u\right\Vert _{m,\Gamma _{1}}^{m}\leq C\left[ \left( \left\Vert
u\right\Vert _{p}^{p}\right) ^{\frac{m\left( 1-s\right) 2}{\left(
2-ms\right) p}}+\left\Vert \nabla u\right\Vert _{2}^{2}\right] .
\end{equation}%
Now, choosing $s$ such that:
$$
0<s\leq \frac{2\left( p-m\right) }{m\left( p-2\right) },
$$
we get:
\begin{equation}\label{choicesm}
\frac{2 m\left( 1-s\right) }{\left( 2-ms\right) p}\leq 1.
\end{equation}
Once the inequality (\ref{choicesm}) is satisfied, we use the classical algebraic inequality:
$$
z^{\nu }\leq \left( z+1\right) \leq \left( 1+\frac{1}{\omega }\right) \left(
z+\omega \right) \;,\quad \forall z\geq 0 \;, \quad  0<\nu \leq 1 \;,\quad \omega \geq 0,
$$
to obtain the following estimate:
\begin{eqnarray}\label{normeLp}
\left( \left\Vert u\right\Vert _{p}^{p}\right) ^{\frac{m\left( 1-s\right) 2}{
\left( 2-ms\right) p}} &\leq & D \left( \left\Vert u\right\Vert
_{p}^{p}+H\left( 0\right) \right)  \nonumber \\
&\leq & D \left( \left\Vert u\right\Vert _{p}^{p}+H\left( t\right) \right)
\;,\quad\forall t\geq 0
\end{eqnarray}%
where we have set $D = 1+ 1/H(0)$.
Inserting the estimate (\ref{normeLp}) into (\ref{estiGamma1}) we obtain the following important 
inequality:
$$
\left\Vert u\right\Vert _{m,\Gamma _{1}}^{m}\leq C\left[ \left\Vert
u\right\Vert _{p}^{p}+\left\Vert \nabla u\right\Vert _{2}^{2}+H\left(
t\right) \right] .
$$
In order to control the term $\left\Vert\nabla u\right\Vert_2^2$ in equation (\ref{derivL3}), we preferely
use (as $H(t) > 0$), the following estimate:
$$
\left\Vert u\right\Vert _{m,\Gamma _{1}}^{m}\leq C\left[ \left\Vert
u\right\Vert _{p}^{p}+\left\Vert \nabla u\right\Vert _{2}^{2}+2 H\left(
t\right) \right] .
$$
which gives finally:
\begin{equation}\label{umgama1}
\left\Vert u\right\Vert _{m,\Gamma _{1}}^{m}\leq C\left[2 d + \left(1+\frac{2}{p}\right) \left\Vert u \right\Vert_p^p -
\left\Vert u_t \right\Vert_2^2 -\left\Vert u_t \right\Vert_{2,\Gamma_1}^2\right] .
\end{equation}
Consequently inserting the inequality (\ref{umgama1}) in the inequality (\ref{derivL3}) we have:
\begin{eqnarray}\label{derivL4}
\frac{dL\left( t\right) }{dt} &\geq &\alpha \left\Vert \nabla u_{t}\right\Vert _{2}^{2} +
\left( r-\frac{\varepsilon \, r \, \left( m-1\right) \, \delta^{-m/\left( m-1\right) } }{m}\right) 
\left\Vert u_{t}\right\Vert_{m,\Gamma _{1}}^{m}  \nonumber \\
&& +  \varepsilon \left( 1 + \frac{r \, C \, \delta ^{m} }{m}\right) \left\Vert u_{t}\right\Vert _{2}^{2} -
\varepsilon \left\Vert \nabla u\right\Vert _{2}^{2}
\\
&&+\varepsilon \left(1-\left(1+\frac{2}{p}\right) \frac{r \, C \, \delta ^{m}  }{m}\right) \left\Vert u\right\Vert _{p}^{p}
+ \varepsilon \left(1+\frac{r \, C \, \delta ^{m}  }{m}\right)\left\Vert u_{t}\right\Vert _{2,\Gamma _{1}}^{2} \nonumber 
\end{eqnarray}
From the inequality (\ref{ineqH}) we have:
$$
- \left\Vert \nabla u\right\Vert_{2}^{2} \geq 2 H(t) - 2 d + \left\Vert u_{t}\right\Vert_{2}^{2}
+ \left\Vert u_{t}\right\Vert _{2,\Gamma_{1}}^{2}
-\frac{2}{p} \left\Vert u\right\Vert _{p}^{p} \; .
$$
Thus inserting it in (\ref{derivL4}), we get the following inequality:
\begin{eqnarray}\label{derivL5}
\frac{dL\left( t\right) }{dt} &\geq &\alpha \left\Vert \nabla u_{t}\right\Vert _{2}^{2} +
\left( r-\frac{\varepsilon \, r \, \left( m-1\right) \, \delta^{-m/\left( m-1\right) } }{m}\right) 
\left\Vert u_{t}\right\Vert_{m,\Gamma _{1}}^{m}  \nonumber \\
&& +  \varepsilon \left( 2 + \frac{r \, C \, \delta ^{m} }{m}\right) \left\Vert u_{t}\right\Vert _{2}^{2} 
+ \varepsilon \left(2+\frac{r \, C \, \delta ^{m}  }{m}\right)\left\Vert u_{t}\right\Vert _{2,\Gamma _{1}}^{2} \\
&& +\varepsilon \left(1-\frac{2 \varepsilon}{p} - \left(1+\frac{2}{p}\right) \frac{r \, C \, \delta ^{m}  }{m}\right) \left\Vert u\right\Vert _{p}^{p}  \nonumber \\
&& + 2 \, \varepsilon \left(H(t) - d \left(1 +\frac{r \, C \, \delta ^{m}  }{m}\right)\right) \nonumber 
\end{eqnarray}
Finally, using the definition of $\alpha_2 \mbox{ and } d$ (see equation (\ref{constant}) and the lemma 
\ref{Vitillaro}), we obtain:  
\begin{eqnarray}\label{derivL6}
\frac{dL\left( t\right) }{dt} &\geq &\alpha \left\Vert \nabla u_{t}\right\Vert _{2}^{2} +
\left( r-\frac{\varepsilon \, r \, \left( m-1\right) \, \delta^{-m/\left( m-1\right) } }{m}\right) 
\left\Vert u_{t}\right\Vert_{m,\Gamma _{1}}^{m}  \nonumber \\
&& +  \varepsilon \left( 2 + \frac{r \, C \, \delta ^{m} }{m}\right) \left\Vert u_{t}\right\Vert _{2}^{2} 
+ \varepsilon \left(2+\frac{r \, C \, \delta ^{m}  }{m}\right)\left\Vert u_{t}\right\Vert _{2,\Gamma _{1}}^{2} \\
&&+\varepsilon \left( \underset{:=c_{0}}{\underbrace{1 -\frac{2}{p} -2 d \left( B\alpha _{2}\right)^{-p}}}-
\left[ \left(1+\frac{2}{p}\right) + 4 d \left( B\alpha _{2}\right)^{-p} \right]\frac{r \, C\, \delta ^{m}}{m}\right) 
\left\Vert u\right\Vert _{p}^{p} \nonumber\\
&&+\varepsilon \left( 2 H (t) +\frac{r \, C\, \delta ^{m}}{m} d)\right) .
\nonumber
\end{eqnarray}%
Setting $\displaystyle c_{0}= 1 -\frac{2}{p}-2d\left( B\alpha_{2}\right) ^{-p} \; $, we have 
$c_0 > 0$ since $\alpha _{2}>B^{-p/(p-2) }$.

We choose now $\delta $ small enough such that:
$$
c_0 - \left[ \left(1+\frac{2}{p}\right) + 4 d \left( B\alpha _{2}\right)^{-p} \right]\frac{r \, C\, \delta ^{m}}{m}>0\quad .
$$
Once $\delta $ is fixed, we choose $\varepsilon $ small enough such that:
$$r-\frac{\varepsilon r\left( m-1\right) }{m}\delta ^{-m/\left(m-1\right) }>0
\text{ and } L(0) > 0 \quad .
$$
Therefore, the inequality (\ref{derivL6}) becomes:
\begin{equation}\label{derivL7}
\frac{dL\left( t\right) }{dt}\geq \varepsilon \eta \left[ H\left( t\right)
+\left\Vert u_{t}\right\Vert _{2}^{2}+\left\Vert u_{t}\right\Vert _{2,\Gamma
_{1}}^{2}+\left\Vert u \right\Vert_{p}^{p} + d \right] \; \mbox{ for some } \eta > 0
\end{equation}
Next, it is clear that, by Young's inequality and Poincar\'{e}'s inequality, we get
\begin{equation}\label{estiL1}
L\left( t\right) \leq \gamma \left[ H\left( t\right) +\left\Vert
u_{t}\right\Vert _{2}^{2}+\left\Vert u_{t}\right\Vert _{2,\Gamma
_{1}}^{2}+\left\Vert \nabla u\right\Vert _{2}^{2}\right]  \; \mbox{ for some } \gamma > 0 .
\end{equation}
Since $H(t) > 0$, we have:
$$
\forall \, t > 0 \,,\, \frac{1}{2} \left \Vert \nabla u\right\Vert_2^2 \leq \frac{1}{p} \left \Vert u\right\Vert_p^p + d \quad .
$$
Thus, the inequality (\ref{estiL1}) becomes:
\begin{equation}\label{estiL1bis}
L\left( t\right) \leq \zeta \left[ H(t) +\left\Vert u_{t}\right\Vert _{2}^{2}+\left\Vert u_{t}\right\Vert _{2,\Gamma_{1}}^{2}+
\left\Vert u \right\Vert_{p}^{p} + d \right]  \;,\; \mbox{ for some } \zeta > 0 .
\end{equation}
From the two inequalities (\ref{derivL7}) and (\ref{estiL1bis}), we finally obtain the differential inequality:
\begin{equation}\label{diffineq}
\frac{dL\left( t\right) }{dt}\geq \mu L\left( t\right)  \;,\; \mbox{ for some } \mu > 0.
\end{equation}%
Integrating the previous differential inequality (\ref{diffineq}) between $0$ and $t$ gives the following estimate for the
function $L$:
\begin{equation}\label{estiL2}
L\left( t\right) \geq L\left( 0\right) e^{\mu t}.
\end{equation}%
On the other hand, from the definition of the function $L$ (and for small values of the parameter $\varepsilon )$, 
it follows that:
\begin{equation}\label{estiLp2}
L\left( t\right) \leq \frac{1}{p}\left\Vert u\right\Vert _{p}^{p}.
\end{equation}%
From the two inequalities (\ref{estiL2}) and (\ref{estiLp2}) we conclude the exponential 
growth of the solution in the $L^{p}$-norm.

\begin{remark} \rm
We recall here that the condition $\displaystyle \int_\Omega u_0(x) u_1(x) dx \geq 0$ appeared in \cite[Theorem 3.12]{GS06}
is unecessary to our result on the exponential growth.
\end{remark}
\vspace*{0.5cm}
{ \bf{Acknowledgment}s}\\
The second author was partially supported by MIRA 2007 project of the R\'egion Rh\^one-Alpes. This author wishes to thank Univ. de Savoie of Chamb\'ery
for its kind hospitality.
Moreover, the two authors wish to thank the referee for his useful remarks and his careful reading of the proofs presented in this paper.

%

\end{document}